\documentclass{ifacconf}

\usepackage{graphicx}      
\usepackage{natbib}        
\usepackage{amssymb}

\usepackage{multicol}
\usepackage{color,enumitem}
\usepackage{bbm,epstopdf, algorithm, multirow}
\usepackage{tabularx}
\usepackage{booktabs}
\usepackage{tabularx}
\usepackage{pifont}
\usepackage{amssymb}
\newcommand{\xmark}{\ding{55}}

\usepackage[noend]{algpseudocode}
\usepackage[normalem]{ulem}
\usepackage{amsmath}
\newcommand{\be}{\begin{equation}}
\newcommand{\ee}{\end{equation}}
\newcommand{\bee}{\begin{eqnarray}}

\ifodd 0
    \newcommand\rev[1]{{\color{magenta}#1}}
    \newcommand{\com}[1]{\textbf{\color{red} (COMMENT: #1)}} 
\else
    \newcommand\rev[1]{{#1}}
    \newcommand{\com}[1]{}
\fi

\newcommand{\eee}{\end{eqnarray}}

\allowdisplaybreaks[4]

\newcommand{\bse}{\begin{subequations}}
\newcommand{\ese}{\end{subequations}}
\newcommand{\defeq}{\overset{\text{\tiny def}}{=}}

\newtheorem{theorem}{Theorem}

\newtheorem{remark}{Remark}

\newcommand\chemfig[1]{{#1}}
\ifodd 0
    
    \newcommand{\comJnl}[1]{\textbf{\color{red} (COMMENT: #1)}} 
\else
    
    \newcommand{\comJnl}[1]{}
\fi

\ifodd 0
    
    \newcommand{\comLcss}[1]{\textbf{\color{red} (COMMENT: #1)}} 
\else
    
    \newcommand{\comLcss}[1]{}
\fi

\ifodd 0 

\else
    
    \newcommand{\accom}[1]{}
\fi

\ifodd 0 
\newcommand{\cnjCom}[1]{{\textbf{\color{red} (CNJ$\to$ #1)}}}
\else
    
    \newcommand{\cnjCom}[1]{}
\fi

\newcommand{\OurAlg}{\textsf{CONFIG}}

\begin{document}
\begin{frontmatter}
\title{\OurAlg: Constrained Efficient Global Optimization for Closed-Loop Control System Optimization with Unmodeled Constraints} 

\author[First,Second]{Wenjie Xu} 
\author[First]{Yuning Jiang} 
\author[Second]{Bratislav Svetozarevic}
\author[First]{Colin N. Jones}
\address[First]{Automatic Control Lab, EPFL, Switzerland (e-mail: wenjie.xu, yuning.jiang colin.jones@epfl.ch)}
\address[Second]{Swiss Federal Laboratories for Materials Science and Technology (EMPA), Switzerland (e-mail: bratislav.svetozarevic@empa.ch)}

\begin{abstract}
In this paper, the \OurAlg{ algorithm}, a simple and provably efficient constrained global optimization algorithm, is applied to optimize the closed-loop control performance of an unknown system with unmodeled constraints. Existing Gaussian process based closed-loop optimization methods, either can only guarantee local convergence~(e.g., \textsf{SafeOPT}), or have no known optimality guarantee~(e.g., constrained expected improvement) at all, whereas the recently introduced \OurAlg{ algorithm} has been proven to enjoy a theoretical global optimality guarantee. In this study, we demonstrate the effectiveness of \OurAlg{ algorithm} in the applications. The algorithm is first applied to an artificial numerical benchmark problem to corroborate its effectiveness. It is then applied to a classical constrained steady-state optimization problem of a continuous stirred-tank reactor. 
Simulation results show that our \OurAlg{ algorithm} can achieve performance competitive with the popular \textsf{CEI}~(Constrained Expected Improvement) algorithm, which has no known optimality guarantee. \rev{As such}, the \OurAlg{ algorithm} offers a new tool, with \rev{both a provable global optimality guarantee and competitive empirical performance}, to optimize the closed-loop control performance for a system with soft unmodeled constraints.
Last, but not least, the open-source code is available as a \textsf{python package} to facilitate future applications.



\end{abstract}


\end{frontmatter}

\section{INTRODUCTION}
\begin{figure}[ht!]
    \centering
    \includegraphics[width=0.5\textwidth]{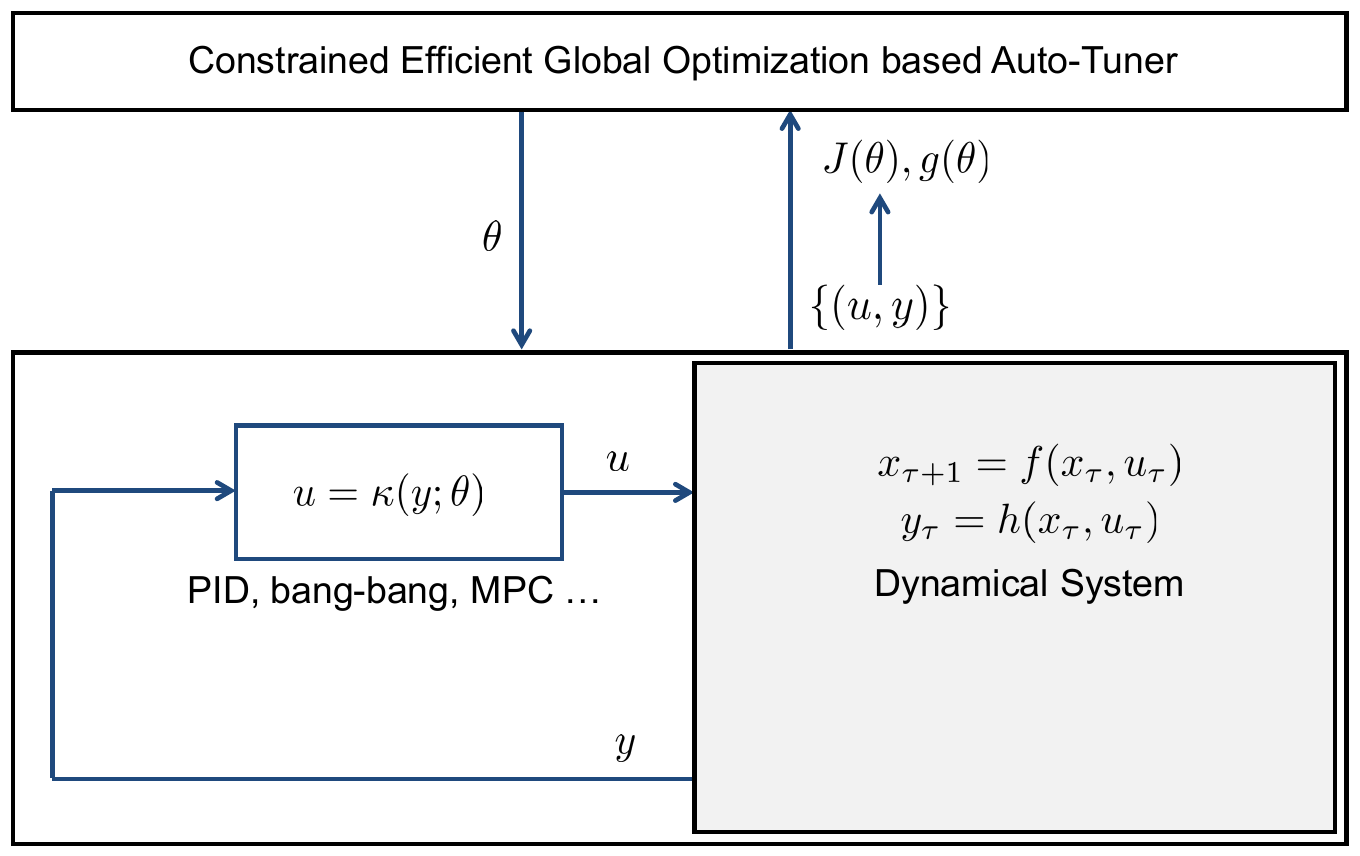}
    \caption{Constrained efficient global optimization based auto-tuner optimizes the performance of a dynamical system controlled by a controller parameterized by $\theta$. By conducting software simulation or physical experiments, we can collect input-ouput trajectories $\{(u,y)\}$ and derive the values of the objective $J(\theta)$ and the constraints $g(\theta)$. On top of the closed-loop dynamical system, we apply our constrained efficient global optimization method to guide the search of optimal feasible parameters $\theta$.}
    \label{fig:ConfigOpt}
\end{figure}

The performance of closed-loop systems can typically be optimized by tuning control parameters~(e.g., gains or set-points) under the guidance of the operational data. Manual tuning of these parameters often takes significant human effort and domain-specific engineering expertise, which may make controller tuning economically infeasible, since the personnel cost may outweigh the benefit from optimized performance. Therefore, algorithms that can automatically adjust these control parameters, without human effort, to optimize closed-loop performance are of great practical interest. 

Closed-loop performance of a control system can typically be defined as functions that take inputs and outputs measured from the closed-loop system as the arguments. For example, in building thermal control, discomfort can be defined as the integration of temperature deviation from a comfort range. While how performance is determined by inputs and outputs is typically clearly understood, the map from control parameters, which are actually tuned, to the closed-loop performance is largely unknown and behaves like a black-box. It is either because an identified system model is unavailable during parameter tuning or the models are so complicated that it is almost impossible to directly construct the explicit mapping from controller parameters to closed-loop performance. 

Therefore, data-driven tuning methods, where the control parameters-to-performance map is considered as a black-box function that can be learned with closed-loop operational data, are receiving more and more research attention. To get closed-loop operation data, we can perform experiments or simulate a model. However, since both experimentation and high-fidelity software simulations are expensive in terms of hardware use or very long computation time. Therefore, it is desired that the tuning algorithms are able to identify a near-optimal set of control parameters with as few experiments/simulations as possible.

Due to the above reasons, efficient global optimization\footnote{Also referred to as Gaussian process optimization~(e.g., in~\cite{srinivas2012information}) or Bayesian optimization~(e.g., in~\cite{frazier2018tutorial}) in machine learning literatures or kriging method~(e.g., in~\cite{jeong2005efficient}) in some engineering literature.}\cite{jones1998efficient} has been widely applied to the closed-loop controller tuning problem, thanks to its black-box modelling capability and superior empirical sample efficiency~\cite{xu2022lower}. Efficient global optimization is a type of sample-efficient derivative-free global optimization method~\cite{jones1998efficient, frazier2018tutorial} that utilizes Gaussian process regression as a surrogate to adaptively search through parameter spaces. In recent work, efficient global optimization has demonstrated potential in controller gain tuning~\cite{lederer2020parameter, duivenvoorden2017constrained,khosravi2019controller, konig2020safety}, MPC tuning~\cite{bansal2017goal, piga2019performance, paulson2020data} and in various other real-world control applications, such as wind energy systems~\cite{baheri2017altitude, baheri2020waypoint}, engines~\cite{pal2020multi} and space cooling~\cite{chakrabarty2021_VCS}.

Meanwhile, in a large variety of applications, another challenge for closed-loop controller tuning is the existence of unmodelled black-box constraints~\cite{xu2022vabo}. For example, when tuning the temperature controller parameters of a building, we need to minimize the energy consumption while keeping the comfort above the occupant-desired level. Both the energy and comfort are unknown black-box functions of the controller parameters. For these applications, we need to extend efficient global optimization based tuning methods to the constrained case.   
 
Fig.~\ref{fig:ConfigOpt} demonstrates the overall structure of the constrained efficient global optimization based controller tuning method, which optimizes the performance of a dynamical system controlled by a controller parameterized by $\theta$. By conducting software simulation or physical experiments, we can collect input-output trajectories $\{(u,y)\}$ and derive the values of the objective $J(\theta)$ and the constraints $g(\theta)$. On top of the closed-loop system, we apply our constrained efficient global optimization method to guide the search of optimal feasible parameters $\theta$. 

To handle these unknown black-box constraints, a variety of efficient global optimization methods with constraints have been developed. We can roughly classify them into different groups based on whether constraint violations are allowed during the optimization process. For the setting where no constraint violations are tolerated, a group of safe Bayesian optimization methods have been developed~\cite{sui2015safe,sui2018stagewise,turchetta2020safe}. However, due to hard safety-critical constraints, these algorithms may get stuck at a local minimum. For another setting where constraint violations are allowed and incur no harm, there exist generic constrained Bayesian optimization methods. For example, a group of popular methods~\cite{gardner2014bayesian,gelbart2014bayesian} encode the constraint information into the acquisition function~(e.g., constrained Expected Improvement). However, there are no known theoretical guarantees on the optimality, constraint violations or convergence rate for this set of methods. 

More recently, a third group of approaches is receiving more  attention~(e.g.,~\cite{xu2022vabo,zhou2022kernelized}), where constraint violations are allowed, but need to be managed well. In this setting, two more recent works adopt a penalty-function approach~\cite{lu2022no} and a primal-dual approach~\cite{zhou2022kernelized} to solve the constrained efficient global optimization problem. The common idea of the two works is the addition of a penalizing term of the constraint violation in the objective and transforming the constrained problem into an unconstrained one. However, both these methods require one to choose values for some critical parameters~(e.g., penalty coefficient~\cite{lu2022no} and dual update step size~\cite{zhou2022kernelized}). The performance of these methods may be impacted heavily by the chosen parameters, \rev{but} no clear guideline on how to select those critical parameters are known 
for implementation. Furthermore, the work of~\cite{lu2022no} only analyzes the penalty-based regret, which is defined as the suboptimality plus penalized constraint violations, but does not derive separate bounds for cumulative regret~(suboptimality) and violations. In practice, improper choice of penalty parameter may lead to the convergence to suboptimal or infeasible solution. Additionally, the cumulative violations considered in~\cite{zhou2022kernelized} are the violations of the cumulative constraint values, not the real total violations. Such a weak form of cumulative violation bounds can not rule out the case of constraint oscillation, where points with severe violations and small constraint values are alternatingly sampled, keeping a low cumulative constraint value. \rev{Besides these works with theoretical guarantees, \cite{priem2020upper} also utilizes upper trust bound to approximate the feasible set and maximize some acquisition function inside this approximate feasible set. However, \cite{priem2020upper} neither gives systematic design of upper trust bound and acquisition function nor provides any theoretical feasibility and optimality guarantees.}

In this paper, we apply the recently introduced \OurAlg{ (\textbf{CON}strained ef\textbf{FI}cient \textbf{G}lobal optimization) algorithm}~\cite{xu2022config}, a simple and provably efficient global optimization algorithm, to optimize the closed-loop control performance of an unknown system with unmodeled constraints. In sharp contrast to existing popoular Gaussian process based closed-loop optimization methods, which either can only guarantee local convergence~(e.g., \textsf{SafeOPT}), or have no known optimality guarantee~(e.g., constrained expected improvement), \OurAlg{ algorithm} has been proven to enjoy a theoretical global optimality guarantee~\cite{xu2022config}. 
Tab.~\ref{tab:comparison} summarizes the comparison of this method to existing state-of-the-art constrained efficient global optimization methods. 
\begin{table*}[htbp!]
    \caption{The comparison with existing constrained efficient global optimization methods.  }
    \centering
     
    \begin{tabular}{|c|c|c|c|c|}
    \hline 
     \textbf{Works} & \begin{tabular}{@{}c@{}}\textbf{Cumulative}\\ \textbf{Regret Bound}\end{tabular}&\begin{tabular}{@{}c@{}}\textbf{Cumulative}\\ \textbf{Violation Bound}\end{tabular} & \begin{tabular}{@{}c@{}}\textbf{Optimality}\\{\textbf{Guarantee}}\end{tabular} & \begin{tabular}{@{}c@{}}\textbf{Infeasibility}\\\textbf{Declaration}\end{tabular}\\
     \hline 
     \begin{tabular}{@{}c@{}}
     \textsf{SafeOPT}\\
     \cite{sui2015safe}, etc.
     \end{tabular}& \xmark & No violation & \begin{tabular}{@{}c@{}}Local convergence\end{tabular}& \xmark \\
     \hline 
     \begin{tabular}{@{}c@{}}
     \textsf{Constrained EI}\\ 
     \cite{gelbart2014bayesian}, etc.
     \end{tabular}& \xmark& \xmark & \xmark &\xmark \\
     \hline
     \begin{tabular}{@{}c@{}}
     \textsf{Primal-Dual method}
     \\
     \cite{zhou2022kernelized}, etc.
    \end{tabular} & \checkmark & weak-form& \xmark & \xmark \\
    \hline 
    \begin{tabular}{@{}c@{}}
      \textsf{Penalty function method}\\
     \cite{lu2022no}, etc.
     \end{tabular}& \multicolumn{2}{c|}{Penalty-based Regret Bound} & \begin{tabular}{@{}c@{}}Global convergence\\ with penalty-dependent rate\end{tabular} & \xmark \\
    \hline 
    \begin{tabular}{@{}c@{}}
        \OurAlg~(Ours)  \\
        \cite{xu2022config} 
    \end{tabular}  & \checkmark& \checkmark & \begin{tabular}{@{}c@{}}Global convergence\\ with rate\end{tabular}  & \checkmark \\
     \hline 
    \end{tabular}
   
    \label{tab:comparison}
\end{table*}

In this study, the algorithm is first applied to artificial numerical problems to indicate its effectiveness. It is then applied to a classical constrained steady-state optimization problem of a continuous stirred-tank reactor.
Simulation results show that \OurAlg{ algorithm} can achieve a performance competitive with the popular \textsf{CEI}~(Constrained Expected Improvement) algorithm, which has no known optimality guarantee. Our algorithm offers a new tool, with a provable global optimality guarantee, to optimize the closed-loop control performance for a system with non-safety critical unmodeled constraints. Last but not least, the open-source code is available as a \textsf{python package} to facilitate future applications.   

The \textbf{contributions} of the paper include:
\begin{enumerate}
\item A \rev{recent} variant of the constrained efficient global optimization method, \OurAlg{, } which achieves provable global optimality, is proposed to do closed-loop control system optimization with unmodeled constraints.

\item The effectiveness and superior performance of the algorithm on artificial problem instances is demonstrated.

\item The method is applied to a constrained steady-state optimization problem of a continuous stirred-tank reactor. 
Simulation results show that this new method achieves performance competitive with the popular state-of-the-art methods that operate without known global optimality guarantees.

\item An open source \textsf{python} implementation of \OurAlg{ is} provided.\footnote{Code is available at: \textsf{https://github.com/JackieXuw/CONFIG}}
\end{enumerate}

\section{Problem Statement}
We consider controlled closed-loop systems of the form
\begin{equation}
x^+ = F(x, \theta),
\label{eq:cl-sys}
\end{equation}
where $x,x^+\in\mathbb{R}^{n_x}$ denote the system state and its update respectively, $\theta\in\Theta\subset \mathbb{R}^{n_\theta}$ are the control parameters (e.g., set-points) to be tuned, and $F(\cdot,\cdot)$ the closed-loop dynamics with some initial state $x_0$.

To determine the system performance, we define a continuous cost function $J(\theta):\mathbb{R}^{n_\theta}\to\mathbb{R}$ to be minimized, which is an unknown/unmodeled black-box function of the parameters $\theta$. We also define $N$ unmodeled constraints on the system outputs that require management during tuning. The $i$-th such constraint is $g_i(\theta):\mathbb{R}^{n_\theta}\to\mathbb{R},i\in[N]$, where the notation $[N]\triangleq\{i\in\mathbb{N}, 1\leq i\leq N\}$.

\begin{remark}
Our formulation can capture two typical application scenarios. 
\begin{itemize}
    \item\textbf{Steady state optimization.} In steady state optimization~(e.g., in~\cite{xu2022vabo}), we assume the state converges to a steady state $x^\infty(\theta)$, which is a function of $\theta$. In this case, the cost function $J(\theta)$ is given as a cost function of the steady state $$J(\theta)\triangleq \ell(x^\infty(\theta)).$$
    and the constraints are also defined accordingly. 
    \item\textbf{Batch performance optimization.} 
   In batch performance optimization, we optimize batch processes over finite time-horizons, say $T_h$. Specifically, we define the objective over a batch trajectory as the integration of a stage cost, 
   \[
   J(\theta)\triangleq\frac{1}{T_h} \int_0^{T_h} l(x(\tau,\theta)) \,\mathrm{d}\tau.
   \] 
   Constraints over the period can also be defined similarly. \rev{For example, in room temperature controller tuning~\cite{fiducioso2019safe}, the aim is to minimize energy consumption, which is the integration of power, subject to constraints including integration of temperature reference tracking error.} 
\end{itemize}
\end{remark}

We formulate the control parameter tuning problem as a black-box constrained optimization problem in Eq.~\eqref{eqn:formulation}. 
\bse\label{eqn:formulation}
\begin{align}
\min_{\theta\in\Theta}  \qquad& J(\theta) \\
\text{subject to:} \qquad & g_i(\theta)\leq 0,\quad \forall i\in[N].\label{eqn:constraint}
\end{align}
\ese

Our \textbf{objective} is to solve the constrained optimization problem~\eqref{eqn:formulation} while managing constraint violations \rev{during the optimization process}. 

\section{Methodology}
In this section, we introduce the preliminaries of using Gaussian process to model the black box functions $J(\theta)$ and $g_i(\theta), i\in[N]$, before introducing the algorithm. 
\subsection{Primer on Gaussian process regression}
As in existing efficient global optimization based controller tuning works~(e.g., \cite{xu2022vabo}), we use Gaussian process surrogate models to learn the unknown functions $J(\theta)$ and $g_i(\theta)$. As in~\cite{chowdhury2017kernelized}, we artificially introduce a Gaussian process $\mathcal{GP}(0, k_i(\cdot, \cdot)), i\in\{0\}\cup[N]$ for the surrogate modelling of the unknown black-box functions $J$ and $g_i$. We also adopt an i.i.d Gaussian zero-mean noise model with noise variance $\lambda$. We use $y_{0,t}$ to denote the measurement of the tuning objective $J$ at step $t$ with the sample $\theta_t$, which is corrupted by independently distributed $\sigma_0$ sub-Gaussian noise. Similarly, we use $y_{i,t}$ to denote the noisy measurement of the tuning constraint $g_i$ at step $t$ with the sample $\theta_t$. We use $\Theta_t$ to denote the sample sequence $(\theta_1,\theta_2,\cdots, \theta_t)$.  
We introduce the following functions of $\theta, \theta^{\prime}$,
\bse
\label{eq:mean_cov}
\begin{align}
\mu_{0,t}(\theta) &=k_{0}(\theta_{1:t}, \theta)^\top\left(K_{0,t}+\lambda I\right)^{-1} y_{0, 1: t}, \\\notag
k_{0,t}\left(\theta, \theta^{\prime}\right) &=k_{0}\left(\theta, \theta^{\prime}\right)\\
&\qquad-k_{0}(\theta_{1:t}, \theta)^\top\left(K_{0,t}+ \lambda I\right)^{-1} k_{0}\left(\theta_{1:t}, \theta^{\prime}\right), \\
\sigma_{0,t}^2(\theta) &=k_{0,t}(\theta, \theta),
\end{align}
\ese
where $k_{0}(\theta_{1:t}, \theta)=[k_0(\theta_1, \theta), k_0(\theta_2, \theta),\cdots, k_0(\theta_t, \theta)]^\top$, $K_{0,t}=(k_0(\theta,\theta^\prime))_{\theta,\theta^\prime\in \Theta_t}$ and $y_{0,1:t}=[y_{0,1}, y_{0, 2},\cdots,y_{0,t}]^\top$. Similarly, we can get $\mu_{i,t}(\cdot), k_{i,t}(\cdot, \cdot), \sigma_{i,t}(\cdot),\forall i\in[N]$ for the constraints. We also introduce the maximum information gain for the objective $f$ as in~\cite{srinivas2012information},
\begin{equation}
    \label{eq:max_inf_gain}
\gamma_{0,t}:=\max_{A \subset \Theta ;|A|=t} \frac{1}{2} \log \left|I+\lambda^{-1}K_{0,A}\right|,
\end{equation}
where $K_{0,A}=(k_0(\theta,\theta^\prime))_{\theta,\theta^\prime\in A}$. 
Similarly, we introduce $\gamma_{i,t},\forall i\in[N]$ for the constraints.

\subsection{\OurAlg{ algorithm for closed-loop controller tuning}}
The \OurAlg{ (\textbf{CON}strained ef\textbf{FI}cient \textbf{G}lobal optimization)} based tuning algorithm is shown in Alg.~\ref{alg:lcb2}, where $\beta^{1/2}_{i,t}$ is a weighting parameter balancing exploration and exploitation. Alg.~\ref{alg:lcb2} adopts the principle of \textit{optimism in the face of uncertainty} for both the objective and constraints. Specifically, we solve a constrained auxiliary problem with the original objective and constraints replaced by their lower confidence bound surrogates, which are posterior mean minus some weighting factor times the posterior standard deviation. Before we solve the auxiliary problem, we check its feasibility. If the auxiliary problem is infeasible, we declare infeasibility for the original controller tuning problem. \cite{xu2022config} gives a way to select $\beta^{1/2}_{i,t}$ so as to guarantee the global optimality for our algorithm.      
\begin{algorithm}[htbp!]
	\caption{\OurAlg{ for Closed-loop Controller Tuning}}
	\begin{algorithmic}[1]
    \For{$t\in[T]$}
      \If
     { $\max_{i\in[N]}\min_{\theta\in\Theta}(\mu_{i,t}(\theta)-\beta^{1/2}_{i,t}\sigma_{i,t}(\theta))>0$}
     \State
     \textbf{Declare infeasibility of~\eqref{eqn:formulation}}. ~\label{alg_line:declare_inf} 
     \EndIf
     \State\label{alg_line:aux_prob} 
     Update controller parameters with 
     \[
     \begin{aligned}
     \theta_t\in{\arg}&\min_{\theta\in \Theta}\quad\mu_{0,t}(\theta)-\beta^{1/2}_{0,t}\sigma_{0,t}(\theta)\\
     \textbf{ subject to }& \mu_{i,t}(\theta)-\beta^{1/2}_{i,t}\sigma_{i,t}(\theta)\leq0,\forall i\in[N].
     \end{aligned}
     \]
     
     \State Run closed-loop system to get measurements of $J$ and $g_i,i\in[N]$.
     \State Update Gaussian process posterior mean and covariance with the new evaluations added. 
    \EndFor
	\end{algorithmic}
\label{alg:lcb2}
\end{algorithm}

We now give the global optimality guarantee of the Alg.~\ref{alg:lcb2}.
\begin{theorem}[\cite{xu2022config}, Theorem 4.2]
Let the same assumptions of Theorem 4.2 in \cite{xu2022config} hold. Then we have, with probability at least $1-\delta,\forall\delta\in(0,1)$, there exists $\tilde{\theta}_T\in\{\theta_1, \theta_2,\cdots, \theta_T\}$, such that,  
\bse
\begin{align}
    J(\tilde{\theta}_T)-J^*&\leq \mathcal{O}\left(\frac{\sum_{i=0}^N\gamma_{i,T}}{\sqrt{T}}\right),\\
[g_i(\tilde{\theta}_T)]^+&\leq \mathcal{O}\left(\frac{\sum_{i=0}^N\gamma_{i,T}}{\sqrt{T}}\right), \quad \forall i\in[N],
\end{align}
\ese
where $J^*$ is the optimal value of the problem~\eqref{eqn:formulation}. 
\label{thm:converge_opt}
\end{theorem}
In the case that the original problem is infeasible, \cite[Theorem 5.1]{xu2022config} also gives a guarantee that infeasibility is declared in finite steps with high probability.


\section{Experiments}
In this section, we will first demonstrate the effectiveness of our algorithm on an artificial numerical problem. We then apply it to a classical constrained steady-state optimization problem of a continuous stirred-tank reactor. 
\subsection{Artificial numerical problem}
We consider the constrained global optimization problem,
\begin{subequations}
\label{eqn:art_num_prob}
\begin{align}
\min_{(\theta^1,\theta^2)\in[-10, 10]^2}\quad   & \cos(2\theta^1)\cos(\theta^2) + \sin(\theta^1),\label{eqn:obj} \\
\text{subject to}\quad &  \cos(\theta^1+\theta^2) - g_\textrm{thr}\leq 0,
\end{align}
\end{subequations}
where $g_\textrm{thr}\in(-1, 1)$. 
To measure the convergence performance of different algorithms, we introduce the concept of constrained regret as 
$$
\min_{t\in[T]}\quad [J(\theta_t)-J^*]^++\sum_{i=1}^N[g_i(\theta_t)]^+,
$$
which is essentially the sum of the positive part of suboptimality and the constraint violations.
\begin{figure}[thbp]
    \centering
    \includegraphics{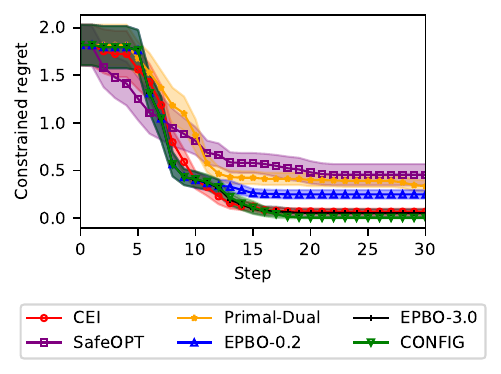}
    \caption{Performance of different algorithms in terms of constrained regret with $g_\textrm{thr}=-0.6$ over 30 instances from randomly selected feasible starting points. Shaded area represents $\pm 0.3$ \textsf{standard deviation}. EPBO-$\rho$ represents EPBO~\cite{lu2022no} method with the penalty coefficient $\rho$.}
    \label{fig:art_ex_convg}
\end{figure}
Fig.~\ref{fig:art_ex_convg} shows the performance of different algorithms in terms of constrained regret, when running on the problem~\eqref{eqn:art_num_prob} with $g_\textrm{thr}=-0.6$ and 30 different feasible starting points. 
Our method can identify the global optimal solution within $20$ steps in this particular example with all the $30$ feasible starting points. 
In sharp contrast, \textsf{SafeOPT} and the primal-dual method can get stuck at a local minimum or infeasible solution, and suffer from a strictly positive \rev{average} constrained regret. The popular \textsf{CEI} method is less efficient and can fail to give an almost optimal solution within $30$ steps. 

As for EPBO method~\cite{lu2022no}, the penalty can be tricky to tune since we do not have information on the gradient of the black-box functions. Randomly selecting a penalty may lead to a suboptimal or infeasible solution if it is not large enough. For example, as shown in Fig.~\ref{fig:art_ex_convg}, EPBO with penalty $0.2$ converges to a strictly positive constrained regret, which means it can fail to find the global optimal feasible solution. On the other hand, too large penalty may also lead to numerical issues. Essentially, \OurAlg{ method} avoids the tuning effort of the penalty parameter and recovers the behavior of EPBO with sufficiently large penalty.   
\subsection{Steady state optimization: Williams-Otto problem}
In this section, we consider the classical Williams-Otto benchmark problem~\cite{del2021real}. In this problem, a continuous stirred-tank reactor~(CSTR) is fed with two pure components \chemfig{A} and \chemfig{B}. The reactor operates at steady state under the temperature $T_\mathrm{r}$. During the reaction, a byproduct \chemfig{G} is also produced. We use $X_\mathrm{A}$ and $X_\mathrm{G}$ to denote the residual mass fractions of \chemfig{A} and \chemfig{G} at the reactor outlet, which need to be managed during the operation. See \cite[Sec.~4.1]{del2021real} for more details. Our goal is to tune the feedrate $F_\mathrm{B}$ of the component \chemfig{B} and the reaction temperature $T_\mathrm{r}$, so as to maximize the economic profit from the reaction. Our problem can then be formulated as, 
\begin{equation}
\begin{aligned}
\min_{F_\mathrm{B}, T_{\mathrm{r}}}\quad& J(F_\mathrm{B}, T_\mathrm{r})\\
\text{subject to}\quad &\textrm{CSTR model~\cite{mendoza2016assessing}}\\
& g_1(F_\mathrm{B}, T_\mathrm{r})\defeq X_\mathrm{A}(F_\mathrm{B}, T_\mathrm{r})-0.12 \leq 0\\
&  g_2(F_\mathrm{B}, T_\mathrm{r})\defeq X_\mathrm{G}(F_\mathrm{B}, T_\mathrm{r})-0.08 \leq 0,\\
& F_\mathrm{B}\in[4, 7], T_\mathrm{r}\in[70, 100]\\ 
\end{aligned}
\end{equation}
where $J(F_\mathrm{B}, T_\mathrm{r})$ is the minimization objective that is opposite to the economic profit, $g_1(F_\mathrm{B}, T_\mathrm{r})$ and $g_2(F_\mathrm{B}, T_\mathrm{r})$ are threshold constraints on the residual mass fractions at the reactor outlet. To measure the quality of a solution $\theta=(F_\mathrm{B}, T_\mathrm{r})$, we use the normalized positive regret plus normalized violations as shown in \begin{equation}
    \frac{[J(\theta)-J^*]^+}{\sigma_J}+\frac{[g_1(\theta)]^+}{\sigma_{g_1}}+\frac{[g_2(\theta)]^+}{\sigma_{g_2}},
\end{equation} 
where $\sigma_J$, $\sigma_{g_1}$ and $\sigma_{g_2}$ are standard deviations of $J, g_1, g_2$ over a set of sampled points. 

\begin{figure}
    \centering
    \includegraphics{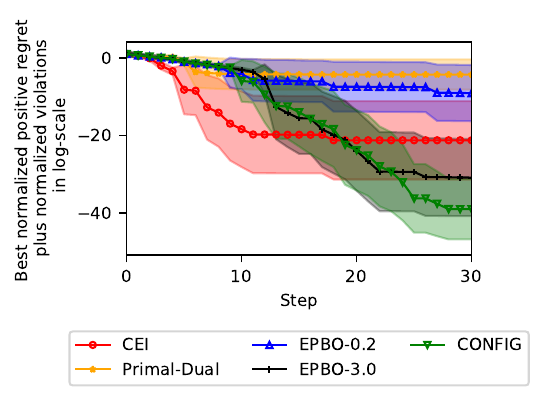}
    \caption{Average performance of different algorithms in log scale from 30 different random starting points. The shaded area represents $\pm 0.5$ \textsf{standard deviation}.}
    \label{fig:WO_funcs_rv_with_epbo_avg}
\end{figure}

Fig.~\ref{fig:WO_funcs_rv_with_epbo_avg} shows the average performance of different algorithms in log scale with 30 different random starting points. We observe that only the \textsf{CEI} method, EPBO-$3.0$ and \OurAlg{ converge} to below $1\times10^{-10}$ within 20 steps, which means they effectively find the approximately constrained optimal solution within 20 steps. Since this example is a black-box, we do not know a feasible solution beforehand and thus, we do not show the result of \textsf{SafeOPT}, which requires an initial feasible solution. As for the EPBO method, when the penalty is set to be too small, it can converge to suboptimal or infeasible solution. Interestingly, with penalty parameter as $3.0$, EPBO performs similarly to \OurAlg{ in} the first 20 steps, but achieves little improvement in the last 10 steps, while \OurAlg{ further} improves the solution.

\section{Computational Toolbox}
To facilitate the future use of \OurAlg{ in} other applications, we build an open-source \textsf{python package}, \textsf{config}, which implements not only the \OurAlg{ algorithm} but also other popular constrained Bayesian optimization methods, to allow more flexibility in the choice of algorithms. 

\begin{figure}[htbp]
    \centering
    \includegraphics[width=\columnwidth]{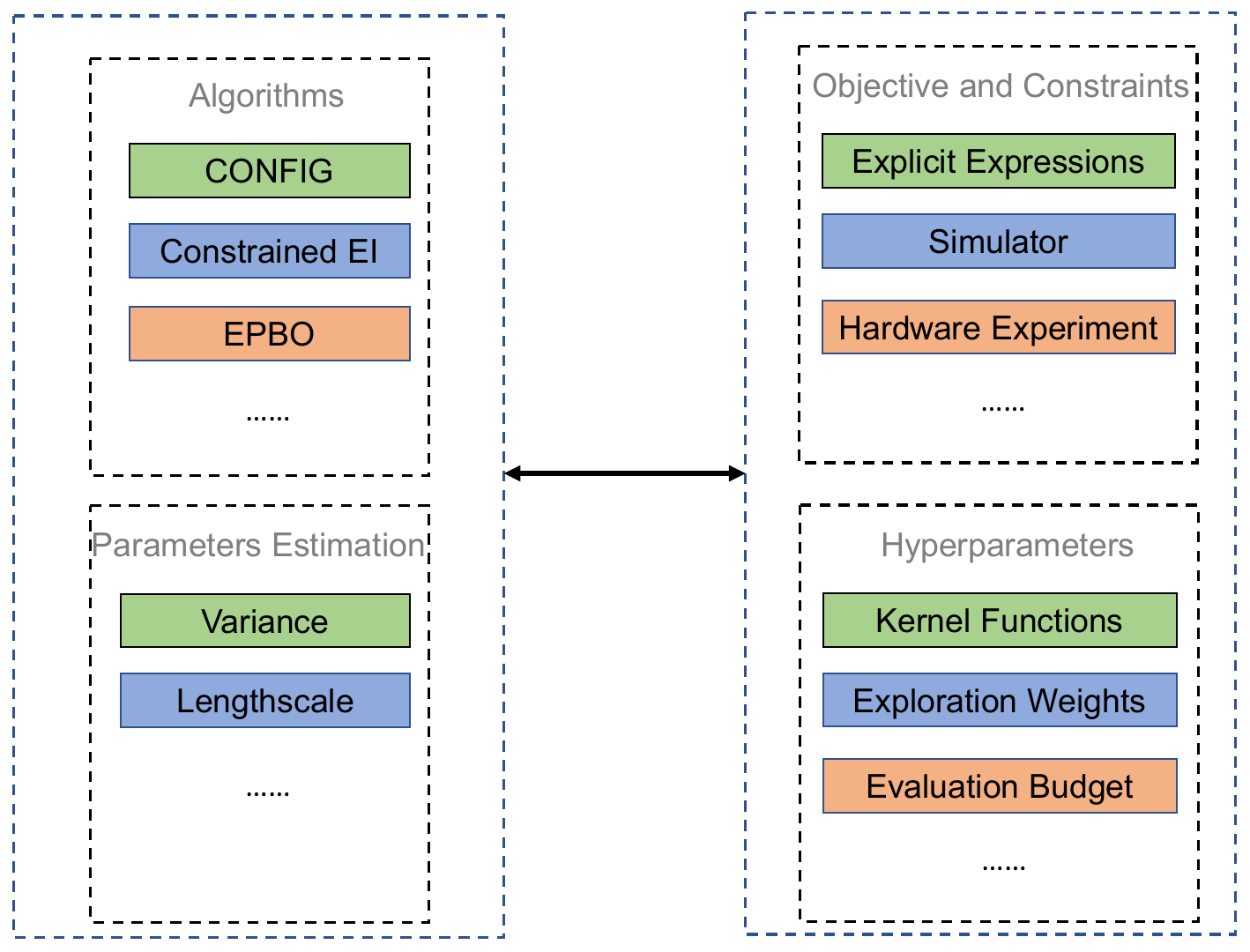}
    \caption{An overview of the major components of our toolbox.}
    \label{fig:comp_toolbox}
\end{figure}
Fig.~\ref{fig:comp_toolbox} shows an overview of the major components for our toolbox to work. Users need to specify the black-box optimization problems, where the objective and constraints can be explicit mathematical expressions, software simulator or hardware experiment. The users also need to specify the choices of some hyperparameters, including kernel functions and evaluation budget. If the user can not provide the parameters in the kernel functions based on domain knowledge, our toolbox also includes simple procedures to estimate the parameters in the kernel functions. The toolbox is extensible, and users can define their own acquisition policies easily using the interface provided.

\section{Conclusion}
In this paper, we have presented the application of the \OurAlg{ algorithm}, a simple and provably efficient constrained global optimization algorithm, to the optimization of the closed-loop control performance of an unknown system with unmodeled constraints. In sharp contrast to currently popular Gaussian process based closed-loop optimization methods, our \OurAlg{ algorithm} enjoys a theoretical global optimality guarantee. Simulation results on an artificial numerical problem and a classical constrained steady-state optimization problem of a continuous stirred-tank reactor
both demonstrate that our \OurAlg{ algorithm} can achieve performance competitive with the popular \textsf{CEI}~(Constrained Expected Improvement) algorithm, which has no known optimality guarantee. Our algorithm offers a new tool, with a provable global optimality guarantee, to optimize the closed-loop control performance for a system with non-safety critical unmodeled constraints. Last but not least, the open-source code is available as a \textsf{python package} to facilitate future applications.   

\begin{ack}
This work was supported by the Swiss National Science Foundation under the NCCR Automation project (grant agreement 51NF40\_180545), the RISK project (Risk Aware Data-Driven Demand Response, grant number 200021\_175627), and by the Swiss Data Science Center (grant number C20-13).
\end{ack}

\bibliography{ifacconf}  

\end{document}